\documentclass[reqno,11pt]{amsart}

\usepackage{amsmath}
\usepackage{amsthm}
\usepackage{mathrsfs}
\usepackage{amsfonts}
\usepackage{parskip} 
\usepackage{hyperref}
\usepackage{epsfig}

\usepackage{amsmath,amssymb,amsfonts,amsthm,amstext,verbatim}
\usepackage{enumerate,color,graphicx,
}
\usepackage{psfrag}  
\usepackage{url}
\usepackage{mdwlist}   
\usepackage{mathrsfs}
\usepackage{hyperref}
\usepackage{epsfig}
\usepackage{etoolbox}

\makeatletter
\patchcmd{\@part}{\Huge}{\normalsize}{}{} 
\makeatother

\numberwithin{equation}{section}

\allowdisplaybreaks

\newtheorem{theorem}{Theorem}[section]

\newtheorem{proposition}[theorem]{Proposition}
\newtheorem{lemma}[theorem]{Lemma}

\newtheorem{definition}[theorem]{Definition}

\newcommand{\ep}{\varepsilon}

\newcommand{\R}{\mathbb{R}}
\newcommand{\N}{\mathbb{N}}
\newcommand{\Z}{\mathbb{Z}}

\newcommand{\dist}{\vert\vert}

\newcommand{\eps}{\ep}

\newcommand{\psap}{\mathsf{P}}

\newcommand{\psaw}{\mathsf{W}}

\newcommand{\esaw}{{\mathbb{E}_{\mathsf{W}}}}

\newcommand{\fpart}{\mathrm{First}}

\newcommand{\lattv}{\Z^d}

\newcommand{\northeast}{\mathrm{NE}}

\newcommand{\vas}{\varsigma}
\newcommand{\vaso}{\varsigma^{\rm loc}}
\newcommand{\gamin}{\Gamma^{{\rm in}}}
\newcommand{\gamout}{\Gamma^{{\rm out}}}
\newcommand{\gaemp}{\ga^{{\rm empty}}}
\newcommand{\ellempty}{l^{{\rm empty}}}

\newcommand{\alphamac}{\alpha'}
\newcommand{\deltamac}{\delta'}
\newcommand{\deltanew}{\varphi}

\newcommand{\ZZ}{\mathbb{Z}}     

\newcommand\hPhi{\mathrm{High}\fpart}

\theoremstyle{remark}

\newcounter{mycount}

\def\mik{1}
\newcommand\cpsfrag[2]{\ifnum\mik=1\psfrag{#1}{#2}\fi}

\newcommand{\saw}{\mathrm{SAW}}

\newcommand{\sawfree}{{\mathrm{SAW}}^{*}}
\newcommand{\sap}{\mathrm{SAP}}

\newcommand{\esm}{\eta}

\newcommand{\thet}{\theta}

\newcommand{\closes}{\text{ closes}}

\newcommand\Ga{\Gamma}     
\newcommand\ga{\gamma}     
\newcommand\Ptmp{P}  
\newcommand\Wtmp{{\rm W}}  
\newcommand\Atmp{{\rm A}}  

\newcommand{\hidden}[1]{}

\newcommand\mcond{\, \middle| \,}
\newcommand\cond{\, | \,}

\newcommand\reverse{\overleftarrow}

\newcommand\pexper{P_{{\rm res}}}
\newcommand\eexper{E_{{\rm res}}}
\newcommand\lmid{\mathrm{l}_{{\rm mid}}}
\newcommand\notick{\mathsf{NoCharm}}

\title[Self-avoiding polygons and walks]{On self-avoiding polygons and walks: \\ the snake method via pattern fluctuation}
\date{}

\author[A.~Hammond]{Alan Hammond}
\address{Departments of Mathematics and Statistics, 
  U.C. Berkeley,
  Berkeley, CA, 94720-3840, U.S.A.}
\email{alanmh@stat.berkeley.edu}
\thanks{The author is supported by NSF grant DMS-$1512908$. 2010 Mathematics Subject Classification. Primary:  60K35.  Secondary: 60D05}

\begin{document}
\maketitle

\vspace{-0.4cm}

\begin{abstract}
For $d \geq 2$ and $n \in \N$, let $\psaw_n$ denote the uniform law on self-avoiding walks of length~$n$ beginning at the origin in the nearest-neighbour integer lattice~$\Z^d$, and write $\Gamma$ for a $\psaw_n$-distributed walk. 
We show that the {\em closing probability} $\psaw_n \big( \dist \Gamma_n \dist = 1 \big)$ that $\Gamma$'s endpoint neighbours the origin is at most $n^{-1/2 + o(1)}$ in any dimension $d \geq 2$.
The method of proof is a reworking of that in~\cite{ontheprob}, which found a closing probability upper bound of $n^{-1/4 + o(1)}$.
A key element of the proof is made explicit and called the {\em snake method}. It is applied to prove the $n^{-1/2 + o(1)}$ upper bound by means a technique of {\em Gaussian pattern fluctuation}.
\end{abstract}
\maketitle

\vspace{1.0cm}

\tableofcontents


\section{Introduction}

Self-avoiding walk was introduced in the 1940s by Flory and Orr \cite{Flory,Orr47} as a model of a long polymer chain in a system of such chains at very low concentration. It is well known among the basic models of discrete statistical mechanics for posing problems that are simple to state but difficult to solve. Two recent surveys are the lecture notes~\cite{BDGS11}  and \cite[Section 3]{Lawler13}. 

\subsection{The model} 

We will denote by $\N$ the set of non-negative integers.
Let $d \geq 2$. For $u \in \R^d$, let $\dist u \dist$ denote the Euclidean norm of $u$. 
A {\em walk} of length $n  \in \N$ with $n > 0$ is a map $\gamma:\{0,\cdots,n \} \to \lattv$ 
such that $\dist \gamma(i + 1) - \gamma(i) \dist = 1$ for each $i \in \{0,\cdots,n-1\}$.
An injective walk is called {\em self-avoiding}. 
A self-avoiding walk~$\gamma$ of length~$n \geq 2$ is said to close (and to be closing)
if $\dist \gamma(n) \dist = 1$. 
When the {\em missing edge} connecting $\gamma(n)$ and $\gamma(0)$
is added, a polygon results.
\begin{definition}
For $n \geq 4$ an even integer, let $\ga:\{0,\ldots,n-1\} \to \Z^d$ be a closing self-avoiding walk. For  $1 \leq i \leq n-1$, let $u_i$ denote the unordered nearest neighbour edge in $\Z^d$ with endpoints $\ga(i-1)$ and $\ga(i)$. Let $u_n$ denote $\ga$'s missing edge, with endpoints $\ga(n-1)$ and~$\ga(0)$. (Note that we have excluded the case~$n=2$ so that $u_n$ is indeed not among the other $u_i$.)
We call the collection of edges 
$\big\{ u_i: 1 \leq i \leq n \big\}$ the polygon of~$\ga$. A self-avoiding polygon in $\Z^d$ is defined to be any polygon of a closing self-avoiding walk in~$\Z^d$. The polygon's length is its cardinality. 
\end{definition}
We will usually omit the adjective self-avoiding in referring to walks and polygons. Recursive and algebraic structure has been used to analyse polygons in such domains as strips, as~\cite{BMB09} describes.

Note that the polygon of a closing walk has length that exceeds the walk's by one. Polygons have even length and closing walks, odd.

Let $\saw_n$ denote the set of self-avoiding walks~$\gamma$ of length~$n$ that start at $0$, i.e., with $\gamma(0) = 0$. 
We denote by $\psaw_n$ the uniform law on $\saw_n$. 
The walk under the law $\psaw_n$ will be denoted by $\Ga$. 
The {\em closing probability} is  $\psaw_n \big( \Gamma \closes \big)$.

In~\cite{ontheprob}, an upper bound on the closing probability of $n^{-1/4 + o(1)}$ was proved in general dimension. 
Without significant modifications, the method used cannot prove an upper bound on this quantity that decays more rapidly than $n^{-1/2 + o(1)}$. In this article, we rework the method in order to reach the conclusion that this latter decay can indeed be achieved.
The next result, which is this conclusion, is the principal result of the present article.


\begin{theorem}\label{t.closingprob}
Let $d \geq 2$. For any $\eps > 0$ and $n \in 2\N + 1$ sufficiently high,
$$
\psaw_n \big( \Gamma \closes \big) \leq n^{-1/2 + \eps} \, .
$$
\end{theorem}


In order to prove  Theorem~\ref{t.closingprob}, we will rework 
 the method of~\cite{ontheprob},  taking the opportunity  to present this method in a general guise, with a view to future applications. Indeed, such an application has already been made, as we explain shortly. It is the introduction and exposition of this general technique which is perhaps the principal advance of the present article. The use of the technique to prove Theorem~\ref{t.closingprob}
 and its relation to the approach of~\cite{ontheprob} is remarked on at the end of Section~\ref{s.indexset}.

The method in question, which is probabilistic in nature, has two elements: the first is a general sufficient condition for proving closing probability upper bounds, and is here explained in a general framework that we call the {\em snake method}. 
In order to apply the snake method to reach an upper bound on the closing probability, such as Theorem~\ref{t.closingprob}, it is necessary to verify the sufficient condition in the method. This step is undertaken in this article, as it was in~\cite{ontheprob}, by a technique of {\em Gaussian pattern fluctuation}. 

\subsection{The snake method via polygon joining}

In~\cite{snakemethodpolygon}, the second application of the snake method appears. The sufficient condition in the method is not verified by pattern fluctuation but instead by a technique of polygon joining. The main result  achieved is now stated.

\begin{theorem}\label{t.thetexist}
Let $d = 2$.  
For any $\eps > 0$, the bound 
$$
\psaw_n \big( \Gamma \closes \big) \leq n^{-4/7 + \eps} 
$$
holds on a set of $n \in 2\N + 1$ of limit supremum density at least $1/{1250}$. 
\end{theorem} 

The decay rate  in the upper bound achieved here   is stronger than that in Theorem~\ref{t.closingprob}, but the result is restricted to dimension $d=2$ and to a positive density set of indices $n \in 2\N + 1$. 

\subsection{An alternative, combinatorial, approach to the closing probability}

Let the walk number $c_n$ equal the cardinality of $\saw_n$. By equation~(1.2.10) 
of~\cite{MS93}, the limit 
$\lim_{n \in \N} c_n^{1/n}$ exists and is positive and finite; it is called the connective constant and denoted by~$\mu$, and we have~$c_n \geq \mu^n$.

Define the polygon number $p_n$ to be the number of length~$n$ polygons up to translation. By (3.2.9) of~\cite{MS93}, 
$\lim_{n \in 2\N} p_n^{1/n} \in (0,\infty)$ exists and equals~$\mu$.

The closing probability may be written in terms of the polygon and walk numbers.
 There are $2n$ closing walks whose polygon is a given polygon of length~$n$, since there are $n$ choices of missing edge and two of orientation. Thus,
\begin{equation}\label{e.closepc}
\psaw_n \big( \Gamma \closes \big) = \frac{2(n+1) p_{n+1}}{c_n} \, ,
\end{equation}
for any $n \in \N$ (but non-trivially only for odd values of $n$).

Since $c_n \geq \mu^n$, an upper bound on $p_n$ of the form $n^{-1-\chi} \mu^n$ implies a closing probability upper bound of $n^{-\chi + o(1)}$. In~\cite{ThetaBound}, such a bound is achieved when $d=2$, $\chi$ is any given value in $(0,1/2)$, and with $n \in 2\N$ in a set of full density in the even integers.  The method used is a  polygon joining technique; as such, two quite different routes to Theorem~\ref{t.closingprob} when $d=2$ (and $n \in \N$ odd is typical) are available.  

\subsection{Conjectural scaling relation and exponent value prediction}

We may hypothesise the existence of exponents $\theta$ and $\xi$
such that $p_n = n^{-\theta + o(1)} \mu^n$ for $n \in 2\N$
and $c_n = n^{\xi + o(1)} \mu^n$ for $n \in \N$. 
Writing  $\psaw_n \big( \Ga \closes\big) = n^{-\psi + o(1)}$, the exponent~$\psi$ would then exist via~(\ref{e.closepc}) and equal $\theta + \xi - 1$.  (We might call $\psi$ the `closing' exponent.)
It is probably fair to say that the existence of the pair $(\theta,\xi)$ is uncontroversial but far from being rigorously established (particularly when $d$ equals two and even more so when it equals three; on results in high dimensions, we comment momentarily). The exponent $\thet$ is predicted to satisfy a relation with the Flory exponent~$\nu$ for mean-squared radius of gyration. The latter exponent is specified by the putative formula $\esaw_n \, \dist \Gamma(n) \dist^2   = n^{2\nu + o(1)}$,
where $\esaw_n$ denotes the expectation associated with~$\psaw_n$ (and where note that $\Gamma(n)$ is the non-origin endpoint of $\Gamma$); in essence, $\dist \Ga(n) \dist$  is supposed to be typically of order $n^{\nu}$. 
The hyperscaling relation that is expected to hold between $\thet$ and $\nu$ is $\thet   =   1 +  d \nu$
where the dimension $d \geq 2$ is arbitrary. In $d = 2$, $\nu = 3/4$ and thus $\theta = 5/2$ is expected. That $\nu = 3/4$ was predicted by the Coulomb gas formalism \cite{Nie82,Nie84} and then by conformal field theory \cite{Dup89,Dup90}. 
We mention also that $\xi = 11/32$ is expected when $d=2$; in light of the $\thet = 5/2$ prediction and~(\ref{e.closepc}),
 $\psaw_n \big( \Ga \closes\big) = n^{-\psi + o(1)}$ with $\psi = 59/32$ is expected. The $11/32$ value was predicted by Nienhuis in \cite{Nie82} and can also be deduced from calculations concerning SLE$_{8/3}$: see \cite[Prediction~5]{LSW}. 

\subsection{Rigorous results in high dimensions}
Hara and Slade~\cite{HS92,HS92b} used the lace expansion to show that $\nu = 1/2$ when $d \geq 5$ by demonstrating that, for some constant $D \in (0,\infty)$,  $\esaw_n \, \dist \Gamma_n \dist^2 - Dn$ is $O(n^{-1/4 + o(1)})$. This value of $\nu$ is anticipated in four dimensions as well, since $\esaw_n \, \dist \Gamma_n \dist^2$ is expected to grow as $n \big( \log n \big)^{1/4}$. (For one article in an extensive recent investigation of Bauerschmidt, Brydges and Slade of the continuous-time weakly self-avoiding walk in $d=4$,   see~\cite{BBSlog}.) Understanding of combinatorial growth rates and closing probability decay is also much more advanced in high dimensions. Indeed, 
in~\cite[Theorem~$1.1(a)$]{HS92},  $c_n$ is shown when $d \geq 5$ to grow as $A \mu^n \big(1 + O(n^{-\eps})\big)$ for some constant $A$ and with $\eps < 1/2$.
The closing probability (and indeed its counterpart for a walk ending at any given displacement from the origin) is shown in~\cite[Theorem~$1.3$]{HS92} to verify $\sum_n n^a \psaw_n (\Ga \closes) < \infty$ for $d \geq 5$ and any $a < d/2 - 1$. This is an averaged version of the stronger assertion that $\psaw_n(\Ga \closes) \leq B n^{-d/2}$ for some large constant $B > 0$: for a proof of this when $d$ is large enough, or on a spread-out lattice, see~\cite[Theorem~$6.1.3$]{MS93}.

\subsection{A suggestion for further reading}
This article has been written in order that it may be read on its own. That said, it does make conceptual sense to view the recent upper bound $p_n \leq n^{-3/2 + o(1)} \mu^n$
(for typical even $n$), and the two snake method applications (yielding Theorems~\ref{t.closingprob}) and~\ref{t.thetexist}), as parts of a whole. The reader who is interested in such a collective presentation should consult the arXiv submission~\cite{CJC}. This article presents in Section~$3$ a heuristic derivation that $\theta \geq 5/2$ when $d=2$ via a technique of polygon joining. Although $\theta = 5/2$ is expected (as we have mentioned), the lower bound derivation provides a useful conceptual framework which is used in~\cite{CJC} to present the various results just mentioned.

\subsection{Structure of the paper}

After some general notation and definitions in Section~\ref{s.four}, we present the general apparatus of the snake method in Section~\ref{s.five} and then use the  method via Gaussian pattern fluctuation in Section~\ref{s.six} to prove Theorem~\ref{t.closingprob}.

\medskip

\noindent{\bf Acknowledgments.} I thank Hugo Duminil-Copin, Alexander Glazman and Ioan Manolescu for stimulating conversations; and a referee for a thorough reading and valuable comments.

\section{Some generalities}\label{s.four}

\subsection{Notation}

\subsubsection{Denoting walk vertices and subpaths}
For $i,j \in \N$ with $i \leq j$, we write $[i,j]$ for  $\big\{ k \in \N:  i \leq k \leq j \big\}$. For a walk $\ga:[0,n] \to \Z^d$ and $j \in [0,n]$, we write $\ga_j$ in place of $\ga(j)$. For $0 \leq i \leq j \leq n$, $\ga_{[i,j]}$ denotes the subpath $\ga_{[i,j]}: [i,j] \to \Z^d$ given by restricting $\ga$. 

\subsubsection{Path reversal}
For $n \in \N$ and a length~$n$ walk $\ga:[0,n] \to \Z^d$, the reversal $\reverse\gamma:[0,n] \to \Z^d$ of $\ga$ is given by $\reverse\ga_j = \ga_{n-j}$ for $j \in [0,n]$.

\subsubsection{Walk length notation} We write $\vert \ga \vert = n$ for the length of any $\ga \in \saw_n$.

\subsubsection{Maximal lexicographical vertex} For any finite $V \subset \Z^d$, we write $\northeast(V)$ for the lexicographically maximal vertex in $V$. The notation is extended to any walk  $\ga:[0,n] \to \Z^d$ by setting $\northeast(\ga)$ equal to $\northeast(V)$ with $V$ equal to the image of~$\gamma$. We will shortly explain the choice of notation~$\northeast$.

\subsubsection{The two-part decomposition}

In the snake method, we represent any given walk $\gamma$
in a {\em two-part} decomposition. This consists of an ordered pair of walks $(\gamma^1,\gamma^2)$ that emanate from
a certain common vertex and that are disjoint except at that vertex.
The two walks are called the {\em first part} and the {\em second part}.
To define the decomposition, consider any walk $\gamma$ of length~$n$. We first mention that the common vertex is chosen to be the lexicographically maximal vertex $\northeast(\gamma)$ on the walk.
Choosing $j \in [0,n]$  
so that $\gamma_j = \northeast(\gamma)$, the walk $\gamma$ begins at $\ga_0$ and approaches $\northeast(\ga)$ along the subwalk $\gamma_{[0,j]}$,
and then continues to its endpoint $\ga_n$ along the subwalk $\gamma_{[j,n]}$. The reversal $\reverse\gamma_{[0,j]}$ of the first walk, and the second walk $\gamma_{[j,n]}$, form a pair of walks that emanate from $\northeast(\gamma)$. (When $j$ equals zero or $n$,
one of the walks is the length zero walk that only visits $\northeast(\gamma)$; in the other cases, each walk has positive length.) The two walks will be the two elements in the two-part decomposition; all that remains is to order them, to decide which is the first part. Associated to each walk is the list of vertices consecutively visited by the walk $\big( \gamma_j,\cdots,\gamma_0 \big)$ and  $\big( \gamma_j,\cdots,\gamma_n \big)$.
These lists may be viewed as elements in $\Z^{d(j+1)}$ and $\Z^{d(n+1-j)}$. 
The first part, $\gamma^1$, is chosen to be the walk in the pair whose list is lexicographically the larger; the second part, $\gamma^2$, is the other walk.

We use square brackets to indicate the two-part decomposition, writing $\ga = [\ga^1,\ga^2]$.

\begin{figure}
    \begin{center}
      \includegraphics[width=0.3\textwidth]{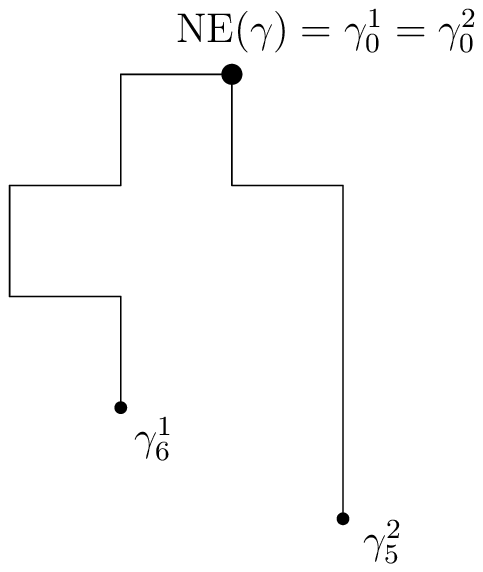}
    \end{center}
    \caption{The two-part decomposition of a walk of length eleven.
}\label{f.firstpart}
  \end{figure}

It is useful to visualise in two dimensions the constructions in our arguments, including this one (see Figure~\ref{f.firstpart}). When $d=2$, we adopt the convention that the lexicographical ordering is specified with the second Euclidean coordinate being recorded before the first. In this way, the lexicographically maximal vertex~$\northeast(\gamma)$ on a walk~$\gamma$ is the vertex that is at least as northerly as any other, and most easterly among those that share its latitude.
The notation $\northeast(\gamma)$ is chosen in light of this 
 `most north then most east' rule. The notation is also used in higher dimensions; we feel permitted to emphasise $d=2$ with the choice of notation because the arguments in this article have almost no dependence on dimension $d \geq 2$.

We write $\big\{ e_1,\cdots,e_d \big\}$
for the consecutive axial Euclidean unit vectors. When $d=2$, $x(u)$ and $y(u)$ will denote coordinates of a point $u \in \Z^2$.

As a small aid to visualization, it is useful to note that if the first $\ga^1$ part of a two-dimensional walk $\ga$ for which $\northeast(\gamma) = 0$ has length~$j \geq 1$, then $\gamma^1:[0,j] \to \Z^2$ satisfies 
\begin{itemize}
\item $\gamma^1_0 = 0$ and $\gamma^1_1 = - e_1$;
\item $y \big( \gamma^1_i \big) \leq 0$ for all $i \in [0,j]$;
\item $\gamma^1_i \not\in \N \times \{ 0 \}$ for any $i \in [1,j]$.
\end{itemize}

\subsubsection{Polygonal invariance}

The following lemma will play an essential role. It is an important indication as to why polygons can be more tractable than walks.

\begin{lemma}\label{l.polyinv}
  For $n \in 2\N + 1$ and $j \in [1,n]$,  let $\chi:[0,n] \to \Z^d$ be a closing walk, and let $\chi'$ be the closing walk obtained from $\chi$ by the cyclic shift $\chi'(i) = \chi\big( j+i  \mod n+1 \big)$, $i \in [0,n]$. Then
  $$
  \psaw_n \Big(\Ga \, \textrm{is a translate of $\chi$} \Big) = 
  \psaw_n \Big(\Ga   \, \textrm{is a translate of $\chi'$} \Big) \, .
  $$
\end{lemma}
\noindent{\bf Proof.} Both sides equal $c_n^{-1}$. \qed
\subsubsection{Notation for walks not beginning at the origin}\label{s.sawfree}
Let $n \in \N$. We write $\sawfree_n$ for the set of self-avoiding walks~$\ga$ of length~$n$ (without stipulating the location~$\ga_0$). We further write $\saw_n^0$ for the subset of $\sawfree_n$ whose elements~$\ga$ have lexicographically maximal  vertex  at the origin, i.e., $\northeast(\ga) = 0$.
 Naturally, an element $\ga \in \saw_n^0$ is said to close (and be closing) if $\dist \ga_n - \ga_0 \dist = 1$. The uniform law on $\saw_n^0$ will be denoted by $\psaw_n^0$. The sets $\saw_n^0$ and $\saw_n$ are in bijection via a clear translation; we will use this bijection implicitly.

\subsubsection{Polygons with northeast vertex at the origin}\label{s.polyconv}

For $n \in 2\N$, let $\sap_n$ denote the set of length $n$ polygons $\phi$ such that $\northeast(\phi) = 0$. The set $\sap_n$ is in bijection with equivalence classes of length~$n$ polygons where polygons are identified if one is a translate of the other. Thus, $p_n =   \vert \sap_n \vert$. 

We write $\psap_n$ for the uniform law on $\sap_n$. 
A polygon  sampled with law $\psap_n$ will be denoted by $\Ga$, as a walk with law $\psaw_n$ is.

There are $2n$ ways of tracing the vertex set of a polygon $\phi$ of length~$n$: $n$ choices of starting point and two of orientation. We now select one of these ways. Abusing notation, and considering $d=2$ for ease of expression but without loss of generality, we may write $\phi$ as a map from $[0,n]$ to $\Z^2$, setting $\phi_0 = \northeast(\phi)$, $\phi_1 = \northeast(\phi) - e_1$, and successively defining $\phi_j$ to be the previously unselected vertex for which $\phi_{j-1}$ and $\phi_j$ form the vertices incident to an edge in $\phi$, with the final choice $\phi_n = \northeast(\phi)$ being made. Note that $\phi_{n-1} = \northeast(\phi) - e_2$.

\subsection{First parts and closing probabilities}
\subsubsection{First part lengths with low closing probability are rare}

\begin{lemma}\label{l.closecard}
Let $n \in 2\N + 1$ be such that, for some $\alphamac > 0$, $\psaw_n \big( \Gamma \closes \big) \geq n^{-\alphamac}$. For any $\deltamac > 0$,  the set of $i \in [0,n]$ for which
$$
 \# \Big\{ \gamma \in \saw^0_n: \vert \gamma^1 \vert = i \Big\}
 \geq n^{\alphamac + \deltamac}  \cdot \# \Big\{ \gamma \in \saw^0_n: \vert \gamma^1 \vert = i  \, , \, \gamma \closes \Big\}
$$
has cardinality at most $2 n^{1- \deltamac}$. 
\end{lemma}
\noindent{\bf Proof.} Note that $\psaw_n \big( \Gamma \closes \big) \geq n^{-\alphamac}$ implies that 
\begin{equation}\label{e.sawcloses}
 \big\vert \saw^0_n \big\vert \leq  n^{\alphamac}  \cdot \# \Big\{ \gamma \in \saw_n: \gamma \closes \Big\} \, .
\end{equation}
Note also that this inequality holds when $\saw_n$ is replaced by $\saw_n^0$.

We have that
$$
\# \Big\{ \gamma \in \saw^0_n: \gamma \closes \Big\} = \sum_{j=0}^n \# \Big\{ \gamma \in \saw^0_n: \gamma \closes \, ,  \, \vert \gamma^1 \vert = j \Big\}
$$
where, by Lemma~\ref{l.polyinv}, each term on the right-hand side has equal cardinality. Writing $Q = Q_{\deltamac} \subseteq [ 0,\ldots,n ]$ for the index set in the lemma's statement, we thus see that
\begin{eqnarray*}
  \big\vert \saw^0_n \big\vert & \geq & \vert Q \vert \cdot
 n^{\alphamac + \deltamac}  \cdot \tfrac{1}{n+1} \, \# \Big\{ \gamma \in \saw^0_n: \gamma \closes \Big\}  \\
 & \geq & \vert Q \vert \cdot  \tfrac{1}{2} \,
 n^{\alphamac -1 + \deltamac}   \, \# \Big\{ \gamma \in \saw^0_n: \gamma \closes \Big\} \, .
\end{eqnarray*}
By~(\ref{e.sawcloses}), or rather by its counterpart for~$\saw_n^0$, we thus find that $\vert Q \vert \cdot \tfrac{1}{2} n^{\alphamac - 1 + \deltamac}$ is at most  $n^{\alphamac}$. \qed

\subsubsection{Possible first parts and their conditional closing probabilities}\label{s.possparts}

Let $d \geq 2$. For $n \in \N$, let $\fpart_n \subseteq \saw_n$ denote  the set of walks  $\gamma:[0,n] \to \Z^d$
whose lexicographically maximal vertex is $\gamma_0 = 0$.
We wish to view $\fpart_n$ as the set of possible first parts of walks $\phi \in \saw^0_m$ of some length $m$ that is at least $n$.
(In the two-dimensional case, we could be more restrictive in specifying $\fpart_n$, stipulating if we wish that any element $\ga$ satisfies $\ga_1 = -e_1$. What matters, however, is only that $\fpart_n$ contains all possible first parts.)

\begin{figure}
    \begin{center}
      \includegraphics[width=0.4\textwidth]{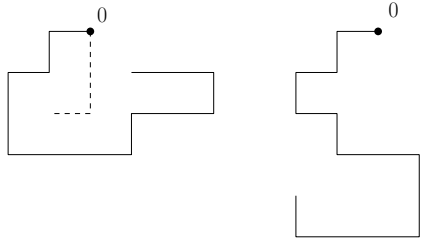}
    \end{center}
    \caption{{\em Left}: the bold $\phi \in \saw_{14}$ and dashed $\ga \in \saw_3$ are such that $[\phi,\ga]$ is a two-part decomposition. Note that $\phi \in \fpart_{14,14+3} \cap \fpart_{14,14+4}^c$. {\em Right}: An element of $\cap_{m = 1}^\infty\fpart_{14,14+m}$.}\label{f.twophi}
  \end{figure}

Note that, as Figure~\ref{f.twophi} illustrates, for given $m > n$, only some elements of $\fpart_n$ appear as such first parts, and we now record notation for the set of such elements (whatever the value of $d \geq 2$). Write $\fpart_{n,m} \subseteq \fpart_n$ for the set of $\gamma \in \fpart_n$ for which there exists an element $\phi \in \saw_{m-n}$ (necessarily with $\northeast(\phi) = 0$) such that $[\gamma,\phi]$ is the two-part decomposition of some element $\chi \in \saw_m$ with $\northeast(\chi) = 0$ (which is to say, $\chi \in \saw_m^0$). 

 In this light, we now define the conditional closing probability 
$$
 q_{n,m}:\fpart_{n,m} \to [0,1] \,   , \, \, \,  \, q_{n,m}(\gamma) = \psaw^0_m \Big( \Gamma \closes \, \Big\vert \, \Gamma^1 = \gamma \Big) \, ,
$$ 
where here $m,n \in \N$ satisfy $m > n$; note that since $\gamma \in \fpart_{n,m}$, the event in the conditioning on the right-hand side occurs for some elements of $\saw_m$, so that the right-hand side is well-defined. 

We also identity a set of first parts with {\em high} conditional closing probability: for $\alpha > 0$, we write
\begin{equation}\label{e.highfirst}
 \hPhi_{n,m}^\alpha = \Big\{ \gamma \in \fpart_{n,m}: q_{n,m}(\gamma) > m^{-\alpha} \Big\} \, .
\end{equation}

\subsection{Classical combinatorial bounds}

\begin{lemma}\label{l.classic}
 Recall that $\mu = \lim_n c_n^{1/n}$.
\begin{enumerate}
\item  There exists a constant $c_{HW} > 0$ such that, for all $n \in \N$, 
$$
  c_n \leq e^{c_{HW} n^{1/2}} \mu^n \, ,
$$
\item There exists a constant $c > 0$ such that, for all $n \in \N$,  $p_n \geq e^{-c n^{1/2}} \mu^n$. 
\end{enumerate}
\end{lemma}
\noindent{\bf Proof: (1).}
This is the 
  Hammersley-Welsh  bound~\cite{hammersleywelsh}, recounted in \cite[Chapter 3]{MS93}.
{\bf (2).}
This bound is \cite[Theorem 3]{kestenone}. 
\qed

\section{The snake method: general elements}\label{s.five}

In this section, we present in a general form  the snake method. 
The method is a  proof-by-contradiction technique which will be applied, alongside other ideas, to prove Theorem~\ref{t.closingprob} in the final Section~\ref{s.six}. The overall strategy of proof in both sections follows that in~\cite{ontheprob}.

The snake method is used to prove upper bounds on the closing probability, and assumes to the contrary that to some degree this probability has slow decay. For the technique to be used, two ingredients are needed.
\begin{enumerate}
\item A charming snake is a walk or polygon $\ga$ many of whose subpaths beginning at $\northeast(\ga)$  have high conditional closing probability, when extended by some common length. It must be shown that charming snakes are not too atypical.
\item A general procedure is then specified in which 
a
charming snake is used to manufacture huge numbers of alternative self-avoiding walks. These alternatives overwhelm the polygons in number and show that the closing probability is very small, contradicting the assumption.
\end{enumerate}

The first step certainly makes use of the assumption of slow decay on the closing probability. It is not however a simple consequence of this assumption. It is this first step that will be carried out in Section~\ref{s.six}, by means of  a technique of Gaussian pattern fluctuation. (In~\cite{snakemethodpolygon}, 
this step is completed via polygon joining to prove Theorem~\ref{t.thetexist}.)

In contrast to these two different approaches that are used to implement the first step, the second step is performed using a general tool, valid in any dimension $d \geq 2$, that we present in this section.  This step draws inspiration from the notion  that reflected walks offer alternatives to closing (or near closing) ones that appears in Madras' derivation~\cite{Madras14} of lower bounds on moments of the endpoint distance under~$\psaw_n$.  

\subsection{The general apparatus of the method}

\subsubsection{Parameters}

The snake method has three exponent parameters:
\begin{itemize}
\item the inverse charm $\alpha > 0$;
\item the snake length $\beta \in (0,1]$;
\item and the charm deficit $\esm \in (0,\beta)$.
\end{itemize}
It has two index parameters:
\begin{itemize}
\item $n \in 2\N + 1$ and $\ell \in \N$, with $\ell \leq n$.
\end{itemize}

\subsubsection{Charming snakes}

Here we define these creatures.

\begin{definition}
Let $\alpha > 0$, $n \in 2\N + 1$,  $\ell \in [0,n]$,
$\gamma \in \fpart_{\ell,n}$, and $k \in [0,\ell]$ with $\ell - k \in 2\N$.  We say that $\ga$ is $(\alpha,n,\ell)$-{\em charming} at (or for) index $k$ if 
\begin{equation}\label{e.charming}
\psaw^0_{k + n - \ell} \Big( \Gamma \closes \, \Big\vert \, \vert \Gamma^1 \vert = k \, , \, \Gamma^1 = \gamma_{[0,k]} \Big) > n^{-\alpha} \, .
\end{equation}
\end{definition}
(The event  that $\vert \Gamma^1 \vert = k$ in the conditioning is redundant and is recorded for emphasis.)
Note that an element $\ga \in \fpart_{\ell,n}$  is $(\alpha,n,\ell)$-charming at  index~$k$ if in selecting a length~$n-\ell$ walk beginning at $0$ uniformly among those choices that form the second part of a walk whose first part is $\ga_{[0,k]}$, the resulting $(k + n - \ell)$-length walk closes with probability exceeding $n^{-\alpha}$. 
(Since we insist that $n$ is odd and that $\ell$ and $k$ share their parity, the length $k + n - \ell$ is odd; the condition displayed above could not possibly be satisfied if this were not the case.) Recalling the definition~(\ref{e.highfirst}),
 note that, for any $\ell \in [0,n]$,  $\ga \in \fpart_{\ell,n}$ is $(\alpha,n,\ell)$-charming at the special choice of index $k = \ell$ precisely when $\ga \in \hPhi_{\ell,n}^{\alpha}$. When $k < \ell$ with $k+n - \ell$ of order~$n$, the condition that 
 $\ga \in \fpart_{\ell,n}$ is $(\alpha,n,\ell)$-charming at  index $k$ is {\em almost} the same as  $\ga_{[0,k]} \in \hPhi_{k,k + n - \ell}^{\alpha}$; (the latter condition would involve replacing $n$ by $k + n -\ell$ in the right-hand side of~(\ref{e.charming})). 

For $n \in 2\N + 1$, $\ell \in [0,n]$, $\alpha,\beta > 0$ and $\esm \in (0,\beta)$, define the charming snake set 
\begin{eqnarray*}
\mathsf{CS}_{\beta,\esm}^{\alpha,\ell,n}  & =  & \Big\{ \ga \in \fpart_{\ell,n}: \textrm{$\ga$ is $(\alpha,n,\ell)$-charming} \\
 & & \quad \textrm{for at least $n^{\beta- \esm}/4$ indices belonging to the interval $\big[ \ell - n^{\beta}, \ell \big]$} \Big\}
 \, .
\end{eqnarray*}
For any element of $\ga \in \fpart_{\ell,n}$, think of an extending snake consisting of $n^{\beta} + 1$ terms
$\big( \gamma_{[0,\ell - n^{\beta}]},\ga_{[0,\ell - n^{\beta} + 1]}, \cdots , \ga_{[0,\ell]} \big)$. If $\ga \in 
\mathsf{CS}_{\beta,\eps}^{\alpha,\ell,n}$, then there are many charming terms in this snake, for each of which there is a high conditional closing probability for extensions of a {\em shared} continuation length~$n - \ell$.

\subsubsection{A general tool for the method's second step}

For the snake method to produce results, we must work with a choice of parameters for which $\beta - \esm - \alpha > 0$. (The method could be reworked to handle equality in some cases.) Here we present the general tool for carrying out the method's second step. This technique was already presented in \cite[Lemma 5.8]{ontheprob}, and our treatment differs only by 
using notation adapted for the snake method with general parameters. 

The tool asserts that, if $\beta - \esm - \alpha > 0$ and even a very modest proportion of snakes are charming, then the closing probability drops precipitously. Recall from Subsection~\ref{s.polyconv}
the notation $\psap_{n+1}$
and a convention for depicting polygons as paths emanating from the northeast vertex.

\begin{theorem}\label{t.snakesecondstep}
Let $d \geq 2$. Set $c = 2^{\tfrac{1}{5(4d + 1)}} > 1$ 
and set $K = 20 (4d+1) \tfrac{\log(4d)}{\log 2}$.
Suppose that the exponent parameters satisfy $\delta = \beta - \esm - \alpha > 0$.
If the index parameter pair $(n,\ell)$ satisfies $n \geq K^{1/\delta}$ and 
\begin{equation}\label{e.afewcs}
 \psap_{n+1} \Big( \Gamma_{[0,\ell]} \in \mathsf{CS}_{\beta,\esm}^{\alpha,\ell,n}  \Big)  \geq c^{- n^{\delta}/2} \, ,
\end{equation}
then
$$
\psaw_n \Big( \Gamma \closes \Big) \leq 2 (n+1) \,  c^{-n^{\delta}/2} \, .
$$ 
\end{theorem}

Note that since the closing probability is predicted to have polynomial decay, the hypothesis~(\ref{e.afewcs}) is never satisfied in practice. For this reason, the snake method will always involve argument by contradiction, with~(\ref{e.afewcs}) being established under a hypothesis that the closing probability is, to some degree, decaying slowly.

\subsection{A charming snake creates huge numbers of reflected walks}
Here is the principal component of the proof of Theorem~\ref{t.snakesecondstep}. 
\begin{proposition}\label{p.avoidance}
Let $d \geq 2$. Set $\delta = \beta - \esm - \alpha$ and suppose that $\delta > 0$. 
Let $\phi \in \mathsf{CS}_{\beta,\esm}^{\alpha,\ell,n}$. With $c  > 1$ and $K > 0$ specified in Theorem~\ref{t.snakesecondstep}, 
we have that, if $n \geq K^{1/\delta}$, then 
$$ 
\# \Big\{ \gamma \in \sawfree_n:  \gamma_{[0,\ell]} = 
 \reverse\phi \Big\} \geq c^{n^{\delta}} \cdot
 \# \Big\{ \gamma \in \sawfree_n:  \northeast(\gamma) = 0 ,  \gamma^1 = \phi \Big\} \, .
$$
\end{proposition}
Note here that walks beginning with the reversal of an element $\phi \in \saw_\ell$ will necessarily not begin at the origin, and thus we employ the notation introduced in Subsection~\ref{s.sawfree}.  

\medskip

\noindent{\bf Proof of Proposition~\ref{p.avoidance}.}
 Let $ \Wtmp$ denote the set of walks $\ga$ of length~$n-\ell$
that originate at $0$ and for which $\northeast(\ga) = 0$. This set is exactly the same as $\fpart_{n-\ell}$ for all $d \geq 2$. We will be using that $\Wtmp$ contains all possible length~$n-\ell$ walks that form the second (rather than the first) part of the two-part decomposition of some walk of at least this length, and thus change notation.

  Let $\Ptmp$ denote the uniform measure on the set $\Wtmp$. We will denote by
  $\Gamma$ a random variable distributed according to $\Ptmp$. 
  In particular, $\Gamma$ is contained in the lower half-space including the origin. (When $d=2$, we mean the region on or below the $x$-axis, and, when $d \geq 3$, the region of non-positive $e_1$-coordinate.)

  We now extend the notion of closing walk by saying that  $\ga'$ {\em closes} $\ga$ if $\ga_0 = \ga'_0$ 
  and the endpoints of  $\ga$ and $\ga'$ are adjacent. 
  We say that $\ga'$ {\em avoids} $\ga$ if no vertex except $\ga_0$ is visited by both $\ga'$ and $\ga$.

We are given $\phi \in \fpart_{\ell,n}$ such that $\phi \in \mathsf{CS}_{\beta,\esm}^{\alpha,\ell,n}$. By definition, we may find indices $j_1 < j_2 < \ldots < j_{\lceil n^{\beta - \esm}/4 \rceil}$ lying in $\big[ \ell - n^{\beta},\ell \big]$ at each of which $\phi$ is $(\alpha,n,\ell)$-charming.

  For $1 \leq i \leq \lceil n^{\beta - \esm}/4 \rceil$, define the events
  \begin{align*}
    A_i & = \Big\{ \Gamma \textrm{ avoids } \, \phi_{[0,j_i]} \Big\} \quad \text{ and } \quad C_i=\Big\{ \Gamma \textrm{ closes } \, \phi_{[0,j_i]}\Big\} \, .
  \end{align*}
  Also, define the set $\Atmp = \big\{ \ga \in \Wtmp : \ga \textrm{ avoids } \phi_{[0,\ell]} \big\}$.
   
  Since $\phi$ is $(\alpha,n,\ell)$-charming at index $j_i$, 
  \begin{align}\label{e.closing.avoiding}
    \Ptmp \Big( \Gamma \closes \, \, \phi_{[0,j_i]} \, \Big\vert \,  \Gamma \text{ avoids } \phi_{[0,j_i]} \Big)
    = \Ptmp \big( C_i \, \big\vert \,  A_i \big) 
    > n^{ - \alpha} \, .
  \end{align}
  Write $k = \lceil  4d \, n^{\alpha} \rceil$. (Note that $k \leq n^{\beta - \esm}/4$ holds for $n$ high enough since $\delta$ is supposed positive.) 
  Any realization $\Ga \in \Wtmp$ is in at most $2d$ events $C_i$.
  Hence, by \eqref{e.closing.avoiding} and the~$A_j$ being decreasing, 
  \begin{align*}
    2d 
    \geq \sum_{i=1}^{k} \Ptmp (C_i)
    \geq \sum_{i=1}^{k} \Ptmp \big(C_i \cond A_i\big) \cdot \Ptmp(A_k)
    \geq 4d \, \Ptmp(A_k) \, .
  \end{align*}
  Therefore, $\Ptmp (A_k) \leq \frac{1}{2}$.
  If the procedure is repeated for indices between $k + 1$ and $2k$, one obtains
  \begin{align*}
    2d
    \geq \sum_{i= k+1}^{2k} \Ptmp (C_i \cond A_k)
    \geq \sum_{i=k+1}^{2k} \Ptmp \big( C_i \cond A_i \big) \cdot \Ptmp \big(A_{2k} \cond A_k \big)
    \geq 4d \, \Ptmp \big( A_{2k} \cond A_k \big) \, ,
  \end{align*}
  and thus $P(A_{2k} \cond A_k ) \leq 1/2$.
  Since $A_{2k}\subset A_k$, we find 
  $$
  \Ptmp (A_{2k}) = \Ptmp (A_k)\Ptmp (A_{2k} \cond A_k) \leq \tfrac14 \, .
  $$
  In these inequalities, we see the powerful bootstrap mechanism at the heart of the snake method, demonstrating that $P(A_{(i+1)k})$ is at most one-half of~$P(A_{ik})$. The mechanism works because
  the method's definitions imply that all walk extensions are of common length~$n-\ell$, and the avoidance conditions are monotone (i.e., the events $A_i$ are decreasing).
  
Indeed,  the procedure may be repeated $ \lfloor \frac{n^{\beta - \esm}}{4k} \rfloor \geq \tfrac{n^{\beta - \esm - \alpha}}{4(4d+1)} - 1$ times. 
Recalling that $\phi = \phi_{[0,\ell]}$, we obtain
  \begin{align*}
    \frac{|\Atmp|}{|\Wtmp|}=\Ptmp\Big(\Gamma\text{ avoids } \phi \Big) \le \Ptmp \big(A_{\lceil n^{\beta - \esm}/4 \rceil} \big) 
    \leq 2^{- \lfloor \tfrac{n^{\beta - \esm}}{4k}\rfloor} 
    \leq 2^{1 - \tfrac{n^{\beta - \esm -  \alpha}}{4(4d+1)}} \, ,
  \end{align*}
whatever the value of $n \in 2\N + 1$.

  \begin{figure}
    \begin{center}
      \includegraphics[width=0.75\textwidth]{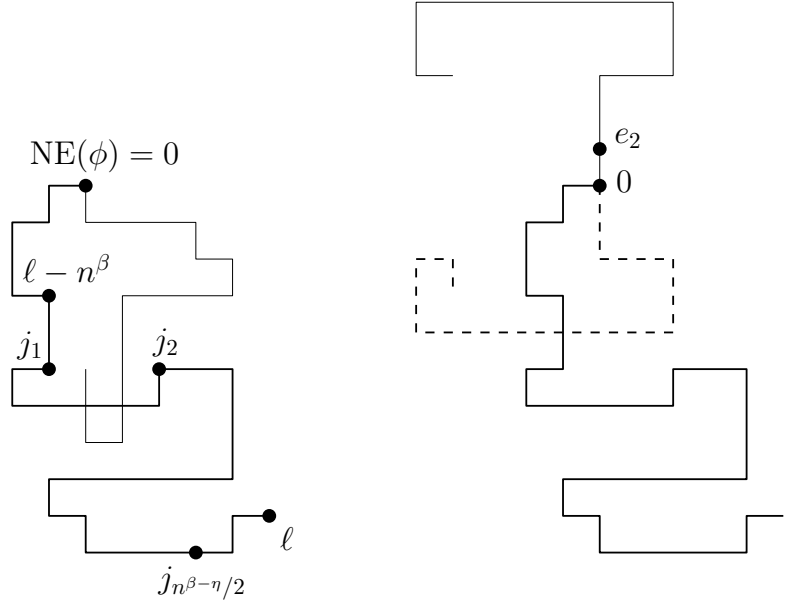}
    \end{center}
    \caption{On the left, $\phi \in \fpart_{\ell,n}$ is depicted in bold. The labels along $\phi$ are indices. A walk in $W \cap A_1 \cap C_1 \cap A_2^c$ is also shown. On the right, another element $\ga \in W$ is depicted with dashed lines. The output of the three-part concatenation procedure used at the end of the proof of Proposition~\ref{p.avoidance} is shown above~$0$.}\label{f.avoided}
  \end{figure}

  The set $\Atmp$ contains all length $n - \ell$ walks $\ga$  
  for which $\big[ \phi , \ga \big]$ is the two-part decomposition of  a walk of length $n$ with $\northeast = 0$. Thus,
$$
 \big\vert \Atmp \big\vert = \# \Big\{ \gamma \in \sawfree_n  :  \northeast(\gamma) = 0 \, , \, \vert \gamma^1 \vert = \ell \, , \, \gamma^1 = \phi \Big\} \, .
$$

  On the other hand, for $\ga \in \Wtmp$, 
  consider the walk obtained by concatenating three paths (and illustrated in Figure~\ref{f.avoided}). When $d=2$, these are: the reversal $\reverse{\phi}$ of $\phi$; 
  the edge $e_2$; and 
the $e_2$-translation of the reflection of $\ga$ in the horizontal axis. When $d \geq 3$, we substitute $e_1$ for $e_2$, and the zero-coordinate hyperplane in the $e_1$-direction for the horizontal axis, to specify these paths.
The walk that results has length $n + 1$ and is self-avoiding. 
  By deleting the last edge of such walks, we obtain at least $|\Wtmp| / 2d$
  walks of length $n$, each of which follows $\reverse\phi$ in its first $\ell$ steps.
  Thus, 
$$
\# \Big\{ \gamma \in \sawfree_n:  \gamma_{[0,\ell]} = \reverse\phi \Big\} \geq |\Wtmp| / 2d \, .
$$
The three preceding displayed equations combine 
prove that the ratio of the cardinalities on the left and right-hand sides in Proposition~\ref{p.avoidance}
is at least $(2d)^{-1} 2^{\tfrac{n^{\beta - \esm -  \alpha}}{4(4d+1)} - 1}$.
The lower bound on $n$ stated in Proposition~\ref{p.avoidance} 
ensures that this last expression is at least $2^{n^{\beta - \eta - \alpha}/{5(4d+1)}}$.
This
completes the proof of this result. \qed

\medskip

\noindent{\bf Proof of Theorem~\ref{t.snakesecondstep}.}
Set $\mathsf{CS}=  \mathsf{CS}_{\beta,\esm}^{\alpha,\ell,n}$. Writing $-\phi_{\ell} + \reverse\phi$
for the translation of  $\reverse\phi$ by $-\phi_{\ell}$, we find that 
\begin{eqnarray*}
 \big\vert \saw_n \big\vert & \geq &   \sum_{\phi \in \mathsf{CS}} \# \Big\{ \gamma \in \saw_n: \gamma_{[0,\ell]} = -\phi_{\ell} + \reverse\phi \Big\} \\
& = & \sum_{\phi \in \mathsf{CS}} \# \Big\{ \gamma \in \sawfree_n: \gamma_{[0,\ell]} = \reverse\phi \Big\} \\
 & \geq &  c^{n^\delta}  \cdot \# \Big\{ \gamma \in \sawfree_n:  \northeast(\gamma) = 0 \, , \, \gamma^1 \in \mathsf{CS} \Big\} \\
 & \geq &  c^{n^\delta}  \cdot \# \Big\{ \gamma \in \sawfree_n:  \northeast(\gamma) = 0 \, , \, \gamma^1 \in \mathsf{CS}  \, , \, \gamma \closes \Big\} \\
 & = &  c^{n^\delta}  \cdot \# \Big\{ \gamma \in \sap_{n+1}:   \gamma_{[0,\ell]} \in \mathsf{CS} \Big\} \geq  c^{n^{\delta}/2}  \cdot    \big\vert  \sap_{n+1} \big\vert \, .  \\
\end{eqnarray*}
Here, we used
 Proposition~\ref{p.avoidance} in the second and~(\ref{e.afewcs}) in the fourth inequalities.
We thus find that
 $p_{n+1}/c_n \leq c^{-n^{\delta}/2}$, and obtain $$\psaw_n \big( \Gamma \closes \big) = 2 (n+1) p_{n+1}/c_n \leq 2 (n+1) c^{-n^{\delta}/2} \, . $$ \qed  

\section{The snake method applied via Gaussian pattern fluctuation}\label{s.six}

Here we prove Theorem~\ref{t.closingprob} by assuming that its conclusion fails and seeking a contradiction. By a relabelling of $\eps > 0$, we may express the premise that the conclusion fails in the form that, for some $\eps > 0$ and infinitely many $n \in 2\N + 1$,
\begin{equation}\label{e.neass}
\psaw_n\big( \Gamma \closes \big) \geq n^{-1/2 + 4\eps} \, .
\end{equation}


We fix the three snake method exponent parameters. Fixing a given choice of $\eps \in (0,1/4)$, we set them equal to  $\alpha = 1/2 - 2\eps$, $\beta = 1/2$ and $\esm  = 0$.
 We will argue that the hypothesis~(\ref{e.afewcs}) is comfortably satisfied, with charming snakes being the~norm.
\begin{proposition}\label{p.snakeg}
There exists a positive constant $C$ such that, if $n \in 2\N + 1$ and $\eps \in (0,1/4)$ satisfy (\ref{e.neass}) as well as $n \geq \max \{ 8^{1/\eps} , C \}$, then there exists  $\ell \in \N$ satisfying $n/4 < \ell < 3n/4$ for which
$$
 \psap_{n+1} \Big( \Gamma_{[0,\ell]} \in \mathsf{CS}_{1/2,0}^{1/2 - 2\eps,\ell,n} \Big) \geq 1 \, - \, 
 n^{-\eps/7}  \, . 
$$
\end{proposition}

\noindent{\bf Proof of Theorem~\ref{t.closingprob}.} 
Given our choice of the three exponent parameters, note that $\delta = \beta - \esm - \alpha$ equals $2\eps$ and is indeed positive. The conclusion of Theorem~\ref{t.snakesecondstep} contradicts~(\ref{e.neass}) if $n$ is high enough. Thus~(\ref{e.neass}) is false for all $n$ sufficiently high. Since $\eps \in (0,1/4)$ may be chosen arbitrarily small, we are done. \qed

\medskip

It remains only to prove Proposition~\ref{p.snakeg}, and the rest of the section is devoted to this proof.
The plan in outline has two steps. In the first, implemented in Lemma~\ref{l.polyclose},
we will infer from the closing probability lower bound~(\ref{e.neass}) a deduction that for typical indices $\ell$ close to $n/2$, the initial (or first part) segment $\Gamma_{[0,\ell]}$
of a $\psap_n$-distributed polygon typically has a not-unusually-high conditional closing probability for a second part extension of length $n - \ell$. This step is a straightforward, Fubini-style statement. In the second step, we work with this deduction to build a charming snake. For this, we want an inference of the same type, but with the second part length remaining fixed, even as the first part length varies over an interval of length close to $n^{1/2}$. Beginning in Subsection~\ref{sec:pattern}, the mechanism of Gaussian fluctuation in local patterns along the polygon is used for this second step. Crudely put, typical deviations in the number of local configurations (type $I$ and type $II$ patterns) in the initial and final ten-percent-length segments of the polygon create an ambiguity in the location marked by the index $\ell \approx n/2$ of order~$n^{1/2}$.
This square-root fuzziness yields charming snakes.  

\subsection{Setting the snake method index parameters}\label{s.indexset}
We now specify the values of $n$ and $\ell$.
The value of $n \in 2\N + 1$ is supposed to satisfy~(\ref{e.neass}) for our given $\eps \in (0,1/4)$,  as well as the bound $n \geq C \vee 8^{1/\eps}$.
Applying Lemma~\ref{l.closecard} with $\alphamac = 1/2 - 4\eps$ and $\deltamac = \eps$, and noting that  $n > 4$ and $n \geq 8^{1/\eps}$ ensures that
$2 n^{1-\eps} < \# \big[ \lceil n/4 \rceil , \lfloor 3n/4  \rfloor \big]$ (since the latter cardinality is at least $n/2 - 1 > 2n^{1-\eps}$), we find that we may select $\ell$ to lie in $\big[ \lceil n/4 \rceil , \lfloor 3n/4 \rfloor  \big]$ and to satisfy 
$$
 \# \Big\{ \gamma \in \saw^0_n: \vert \gamma^1 \vert = \ell \Big\}
 < n^{1/2 - 3\eps}  \cdot \# \Big\{ \gamma \in \saw^0_n: \vert \gamma^1 \vert = \ell  , \gamma \closes \Big\} \, ,
$$
or equivalently
$$
\psaw^0_n \Big(  \Gamma \closes \, \Big\vert \, \vert \Gamma^1 \vert = \ell   \Big) > n^{-1/2 + 3\eps} \, .
$$
The value of $\ell$ is so fixed henceforth in the proof of Proposition~\ref{p.snakeg}.

\begin{lemma}\label{l.polyclose}
Recalling~(\ref{e.highfirst}), we have that
$$
 \psap_{n+1} \Big( \Gamma_{[0,\ell]} \not\in \hPhi_{\ell,n}^{1/2 - 2\eps} \Big) \leq n^{-\eps} \, .
$$
\end{lemma}
\noindent{\bf Proof.} 
Note that $\Gamma_{[0,\ell]}$ under $\psap_{n+1}$ shares its law with the first part $\Gamma^1$ under $\psaw^0_n \big( \cdot \, \big\vert \, \Gamma \closes \, , \, \vert \Gamma^1 \vert = \ell \big)$. For this reason, the statement may be reformulated
\begin{equation}\label{e.reform}
 \psaw^0_n \Big(  \Gamma^1 \not\in \hPhi_{\ell,n}^{1/2 - 2\eps} \, \Big\vert \, \vert \Gamma^1 \vert = \ell \, , \, \Gamma \closes \Big) \leq n^{-\eps} \, .
\end{equation}
To derive (\ref{e.reform}), set $p$ equal to its left-hand side. Note that
\begin{eqnarray*}
  & & \# \Big\{ \gamma \in \saw^0_n: \vert \gamma^1 \vert = \ell   \Big\} \\
 & \geq &   \# \Big\{ \gamma \in \saw^0_n: \vert \gamma^1 \vert = \ell \, , \, \gamma^1 \not\in \hPhi_{\ell,n}^{1/2 - 2\eps}   \Big\} \\
 & = & \sum_{\phi \in \fpart_{\ell,n} \setminus \hPhi_{\ell,n}^{1/2 -  2\eps}}  \# \Big\{ \gamma \in \saw^0_n:  \gamma^1  = \phi   \Big\} \\
 & \geq & \sum_{\phi \in \fpart_{\ell,n} \setminus \hPhi_{\ell,n}^{1/2 - 2\eps}} n^{1/2 - 2\eps} \cdot  \# \Big\{ \gamma \in \saw^0_n:  \gamma^1  = \phi \, , \, \gamma \closes  \Big\} \\
 & = &   n^{1/2 -  2\eps} \cdot  \# \Big\{ \gamma \in \saw^0_n: \vert \gamma^1 \vert = \ell \, , \, \gamma^1 \not\in \hPhi_{\ell,n}^{1/2 - 2\eps} \, , \, \ga \closes  \Big\} \\ 
& = &   n^{1/2 - 2\eps} \cdot  p \cdot \# \Big\{ \gamma \in \saw^0_n: \vert \gamma^1 \vert = \ell  \, , \, \ga \closes  \Big\} \, .
\end{eqnarray*}
The second inequality exploits the concerned $\phi$  not lying in $\hPhi_{\ell,n}^{1/2 -  2\eps}$. We learn that
$$
 n^{1/2 - 2\eps} \cdot  p \leq \psaw^0_n \Big( \Ga \closes \, \Big\vert \, \vert \Gamma^1 \vert = \ell  \Big)^{-1} \, , 
$$
whose right-hand side we know to be at most $n^{1/2 - 3\eps}$. Thus, $p \leq n^{-\eps}$, and we have verified~(\ref{e.reform}). \qed

\medskip

\noindent{\em Remark.}
Lemma~\ref{l.polyclose} may be compared to \cite[Lemma 5.5]{ontheprob} with $k=0$.
The former result states that $\hPhi_{\ell,n}^{1/2 - 2\eps}$ membership by~$\Gamma_{[0,\ell]}$ is the norm under~$\psap_{n+1}$, while the latter merely asserts that a comparable membership is not unlikely (having probability at least $n^{-1/4 + 2\eps}$).
It may be possible to improve the inequality in~\cite[(5.7)]{ontheprob}
to reflect the fact that the conditioning on closing leading from the law $\psaw_n$
to $\psap_{n+1}$ reweights the measure on first parts proportionally in accordance
with the conditional closing probability given the first part. Analysing this reweighting may lead to a replacement of the right-hand side of \cite[(5.5)]{ontheprob} by a term of the form $1 - n^{-o(1)}$ and, alongside other suitable changes, permit a derivation of Theorem~\ref{t.closingprob}.

\subsection{Patterns and shells}\label{sec:pattern}

Patterns are local configurations in self-avoiding walks that are the subject of a famous theorem~\cite{kestenone} due to Kesten that we will shortly state. For our present purpose, we identify two particular patterns.

\begin{figure}
  \begin{center}
    \includegraphics[width=0.4\textwidth]{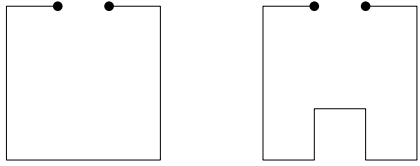}
  \end{center}
  \caption{An example of type I and II patterns for $d =2$.}
  \label{f.patternstwo}
\end{figure}

\begin{definition}[Type I/II patterns] 
  A pair of type I and II patterns is a pair of self-avoiding walks $\chi^I$, $\chi^{II}$,  
  both contained in the cube $[0,3]^d$, with the properties that
  \begin{itemize*}
  \item $\chi^I$ and $\chi^{II}$ both visit all vertices of the boundary of $[0,3]^d$,
  \item $\chi^I$ and $\chi^{II}$ both start at   $\big( 1,3,1, \cdots , 1 \big)$ and end at $\big( 2,3, 1, \cdots, 1 \big)$,
  \item the length of $\chi^{II}$ exceeds that of $\chi^{I}$ by two.
  \end{itemize*}
\end{definition}
Figure \ref{f.patternstwo} depicts examples of such patterns for $d =2$. 
The existence of such pairs of walks for any dimension $d \geq 2$ may be easily checked, and no details are given here. 
Fix a pair of type $I$ and $II$ patterns henceforth. 

A pattern $\chi$ is said to occur at step $k$ of a walk $\ga$ 
if $\ga_{[k, k + |\chi|]}$ is a translate of $\chi$ (where recall that $\vert \chi \vert$ is the length of~$\chi$).
A {\em slot} of $\ga$ is any translate of $[0,3]^d$ containing $\ga_{[k, k + |\chi|]}$ where a pattern~$\chi$ of type $I$ or $II$ occurs at step $k$ of $\gamma$.
Note that the slots of $\ga$ are pairwise disjoint. 

In~\cite{ontheprob}, the notion of shell was introduced. A shell is an equivalence class of self-avoiding walks under the relation that two walks are identified if one may be obtained from the other by changing some patterns of type~$I$ to be of type~$II$ and {\em vice versa}. The walks in a given shell share a common set of slots, but they are of varying lengths. The shell of a given walk $\ga$ is denoted~$\vas(\ga)$.

Consider two walks $\ga^1$ and $\ga^2$ with $\ga^1_0 = \ga^2_0$. Recall that $\ga^2$ {\em avoids} $\ga^1$ if  no other vertex is visited by both walks. 
The next fact is crucial to our reasons for considering shells; its almost trivial proof is omitted.
\begin{lemma}\label{l.sharedavoidance}
For some $m \in \N$, let $\ga \in \saw_m$ and let  $\ga' \in \vas(\ga)$.
A walk beginning at $\ga_0$ avoids $\ga$ if and only if it avoids $\ga'$. 
\end{lemma}

\begin{figure}
  \begin{center}
    \includegraphics[width=1.0\textwidth]{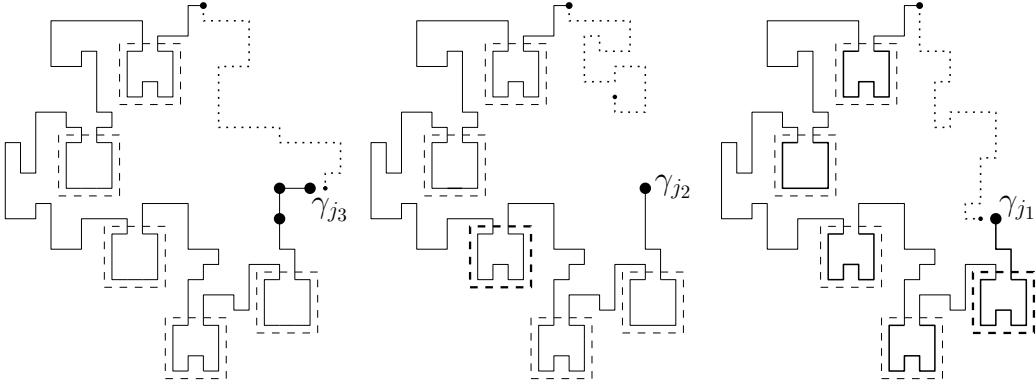}
  \end{center}
  \caption{In this figure, we explain in outline the method. Given $n \in 2\N + 1$, the index $\ell$ has been fixed between $n/4$ and $3n/4$ so that the bound~(\ref{e.reform}) holds. As we have seen, the vast majority of indices in this range satisfy this bound. This means that when we draw a  length~$n+1$ polygon $\ga$ and mark with a black dot each vertex $\ga_j$, $n/4 \leq j \leq 3n/4$, with the property that 
  $\psaw^0_n \big( \Ga \closes \, \big\vert \, \vert \Ga^1 \vert = j , \Ga^1 = \ga_{[0,j]} \big) \geq n^{-1/2 + 2\eps}$,  
   most such~$\ga$ will appear with black dots in most of the available spots.  The left-hand sketch represents such a typical~$\ga$ and three of its many black dots. The pointed second part shows a sample of the law $\psaw^0_n \big( \cdot \, \big\vert \, \vert \Ga^1 \vert = j_k , \Ga^1 = \ga_{[0,j_k]} \big)$, with $k=3$, one that happens to close~$\ga_{[0,j_3]}$. The second part being sampled has length $n - j_3$. If we sample instead this law with $k=2$, then the second part has a greater length, $n - j_2$, which equals $n-j_3 + 2$ in this instance. However, to construct a charming snake, we want this length to stay the same as we move from one black dot to the next. Pattern exchange is the mechanism that achieves this. By turning one type~$I$ pattern into a pattern of type~$II$,  we push two units of length into the first part, so that, in the middle sketch, the random second part has the original length $n - j_3$. The process is iterated in the right-hand sketch. The first part is akin to a belayer who takes in rope, storing it in accumulating type~$II$ patterns, so that the second part climber maintains a constant length of rope. This process of pattern exchange can be maintained for an order of $n^{1/2}$ steps, because the Gaussian fluctuation between the two types of pattern means that the process of artificially altering pattern type does not push the system out of its rough equilibrium when the number of changes is of this order. In this way, black dots also mark charming snake terms for a snake of length of order~$n^{1/2}$.}\label{f.snakemethodpattern}
\end{figure}

The reader may now wish to view Figure~\ref{f.snakemethodpattern} and its caption for an expository overview of the snake method via Gaussian pattern fluctuation. We mention also that this Gaussian fluctuation has been utilized in~\cite{JROSTW} to prove a $n^{1/2 - o(1)}$ lower bound on the absolute value of the writhe of a typical length~$n$ polygon.

We will make some use of the notion of shell, but will predominantly consider 
a slightly different definition, the $(n+1)$-local shell, which we now develop.
This new notion concerns polygons rather than walks. 

Recall that the parameter $n \in 2\N +1$ has been fixed in Subsection~\ref{s.indexset}. (The upcoming  definitions do not require that $n \in 2\N + 1$ be fixed in this particular way in order to make sense, but, when we make use of the definitions, it will be for this choice of $n$.)
Define an equivalence relation $\sim$ on $\sap_{n+1}$ as follows. For any $\gamma \in \sap_{n+1}$, let $\gaemp$ denote the polygon in the shell of $\gamma$ that has no type $II$ patterns, (formed by switching every type $II$ pattern of $\ga$ into a type $I$ pattern). Thus, $\gaemp \in \sap_{n + 1 - 2T_{II}(\ga)}$, where $T_{II}(\ga)$ denotes the total number of type~$II$ patterns in $\ga$. A type $II$ pattern contains thirteen edges (in $d=2$; at least this number in higher dimensions), and these patterns are disjoint, so $T_{II}(\ga) \leq (n+1)/13$. Thus, the length of $\gaemp$ is at least $11(n+1)/13$. Let $S_1$ denote the set of slots in $\gamma$ that are slots in $\gaemp_{[0,(n+1)/10]}$, and, writing $\ellempty$ for the length of $\gaemp$, let $S_2$ denote the set of slots in $\gamma$ that are slots in $\gaemp_{[\ellempty - (n+1)/10,\ellempty]}$. Note that $S_1$ and $S_2$ are disjoint. We further write $N_I(\ga)$ and $N_{II}(\ga)$ for the number of patterns of the given type in the slots $S_1 \cup S_2$, and $N_{I}^1(\ga)$ and $N_{I}^2(\ga)$ for the number of type~$I$ patterns occupying slots in $S_1$ and in $S_2$; and similarly for $N_{II}^1(\ga)$ and $N_{II}^2(\ga)$. 

For $\ga,\ga' \in \sap_{n+1}$, we say that $\ga \sim \ga'$ 
if $\gamma'$ may be obtained from $\gamma$ by relocating the type $II$ patterns of $\gamma$ contained in the set of slots $S_1 \cup S_2$ for $\gamma$ to another set of locations among these slots. The relation $\sim$ is an equivalence relation, because the polygon $\gaemp$ formed by filling all the slots of $\ga$ with type $I$ patterns is shared by related polygons, so that the value of $S_1 \cup S_2$ is equal for such polygons. 
Elements of a given equivalence class have common values of length, $N_I$ and $N_{II}$, but not of $N_I^i$ or $N_{II}^i$ for $i \in \{1,2\}$. 
We call the equivalence classes $(n+1)$-local shells: the parameter $n+1$ appears to denote the common length of the member polygons, and the term {\em local} is included to indicate that members of a given class may differ only in locations that are close to the origin (in the chemical distance, along the polygon). Complementing the notation $\varsigma(\ga)$ for the shell of $\ga$,
 write $\vaso(\ga)$ for the $(n+1)$-local shell of $\ga \in \sap_{n+1}$.

For $\deltanew > 0$, write $\mathcal{G}_{n+1,\deltanew}$ for the set of $(n+1)$-local shells $\sigma \subseteq  \sap_{n+1}$ such that each of the quantities $\vert S_1(\sigma) \vert$, $\vert S_2(\sigma) \vert$, $N_I$ and $N_{II}$ is at least $\deltanew (n+1)$. Such ``good'' shells are highly typical if $\deltanew > 0$ is small, as we now see.

\begin{lemma}\label{l.kptconseq}
There exist constants $c > 0$ and $\deltanew > 0$ such that 
$$
\psap_{n+1} \Big( \vaso(\Gamma) \in \mathcal{G}_{n+1,\deltanew} \Big) \, \geq \, 1 - e^{-cn} \, .
$$ 
\end{lemma}
\noindent{\bf Proof.}
By Kesten's pattern theorem \cite[Theorem 1]{kestenone}, 
there exist constants $c > 0$ and $\deltanew >0$ such that,
for any odd $n \geq d 3^d$,
\begin{equation}\label{e.patterns}
  \psaw_{n+1} \Big( T_I(\Ga) \leq \deltanew n \Big) 
  \leq 
  e^{-c n}  \, ,
\end{equation}
where naturally $T_I(\Ga)$ denotes the number of type~$I$ patterns in $\Gamma$. 

Note that every slot of $\Gamma_{[0,(n+1)/10]}$
is also a slot in $S_1(\Gamma)$: after all, the distance to be traversed to reach a given vertex in 
$\Gamma^{{\rm empty}}$ is at least as great along $\Gamma$ as it is along  $\Gamma^{{\rm empty}}$, because of the removal of type~$II$ patterns by which  
$\Gamma^{{\rm empty}}$  is formed from $\Gamma$.
Thus, 
\begin{eqnarray*}
 & & \psap_{n+1} \Big( S_1(\Ga) \text{ contains fewer than $\deltanew (n+1)$ type I patterns} \Big) \\
  & \leq & \psap_{n+1} \Big( \Ga_{[0, (n+1)/10]} \text{ contains fewer than $\deltanew (n+1)$ type I patterns} \Big) \,  .
\end{eqnarray*}

Set
$$
 \omega =  \psaw_{(n+1)/10} \Big( \Ga\text{ contains fewer than $\deltanew (n+1)$ type I patterns} \Big) \, ,
$$
so that~(\ref{e.patterns}) implies $\omega \leq e^{-c(n+1)/10}$.
  There exist constants $\deltanew ,c > 0$ such that 
  \begin{eqnarray*}
    & & \psap_{n + 1} \Big( \Ga_{[0,(n + 1)/10]} \text{ contains fewer than $\deltanew (n + 1)$ type I patterns} \Big) \nonumber \\
    &\le &   \frac{c_{(n+1)/10} c_{9(n+1)/10}}{p_{n+1}} \cdot \omega 
   \, \le \, \exp \big\{ (4 (10)^{-1/2} c_{{\rm HW}} + c) (n+1)^{1/2} \big\} \cdot \omega \nonumber
    \, \leq \, e^{-c n},
  \end{eqnarray*}
where the second inequality relies on both parts of Lemma~\ref{l.classic}, and the third entails a relabelling of $c > 0$. Thus, $\psap_{n+1} \big( N_{I}^1 \leq \deltanew (n+1) \big) \leq e^{-cn}$. The same holds for the quantity $N_{II}^1$. Considering 
$\Ga_{[9(n+1)/10, n+1]}$ in place of $\Ga_{[0,(n + 1)/10]}$, the same conclusion may be reached about $N_I^2$ and $N_{II}^2$. 
  It follows that 
  \begin{align*}
    \psap_{n+1} \Big( \min \big\{ \vert S_1 \vert , \vert S_2 \vert , N_I , N_{II} \big\} < \deltanew (n+1) \Big) < 4e^{-c n} \, .
  \end{align*}
This completes the proof. \qed  

\medskip

In the next lemma, we see how, for any $\ga \in \mathcal{G}_{n+1,\deltanew}$, the mixing of patterns that occurs when an element of $\vaso(\ga)$ is realized involves an asymptotically Gaussian fluctuation in the pattern number $N_I^1$. The statement and proof are minor variations of those of \cite[Lemma 3.5]{ontheprob}.

\begin{lemma}\label{l.balanced}

For any $\deltanew > 0$, there exists $c > 0$ and $N \in \N$ such that, for $n \geq N$ and $\ga \in \mathcal{G}_{n+1,\deltanew}$,

\medskip

\noindent{\bf (1)}  if $k \in \N$ satisfies 
  $\left| k - \frac{T_I |S_1|}{|S_1|+|S_2|}\right|  \leq  n^{1/2} (\log n)^{1/4}$,  then
$$ 
\psap_{n+1}\Big( N^1_I(\Ga) = k \, \Big\vert \, \Gamma \in \vaso(\ga) \Big)
    \geq n^{-1/2} \exp \big\{ - c (\log n)^{1/2} \big\} \, ;    
$$
\noindent{\bf (2)} and, for any  $g \in [1/8,1/2]$,
$$ 
\psap_{n+1}\bigg( \left| N^1_I(\Ga) -  \frac{N_I |S_1|}{|S_1|+|S_2|}\right|  \geq  n^{1/2} (\log n)^{g} \, \bigg\vert \, \Gamma \in \vaso(\ga) \bigg)
    \leq  \exp \big\{ - c (\log n)^{2g} \big\} \, .    
$$
\end{lemma}

\noindent{\bf Proof.}
Fix a choice of $\ga$ in $\mathcal{G}_{n+1,\deltanew}$.   If $\Ga$ is distributed according to $\psap_{n+1}$ conditionally on $\Gamma \in \vaso(\ga)$, then 
  $N_I$ type $I$ patterns and $N_{II}$ type $II$ patterns are distributed uniformly in the slots of $S_1 \cup S_2$. 
  Thus, for $k \in \big\{0, \cdots, \vert S_1 \vert \big\}$,
  \begin{align}\label{e.Tdistrib}
    \psap_{n+1}\Big( N^1_I(\Ga) = k  \, \Big\vert \, \Gamma \in \vaso(\ga) \Big) 
    = \frac{\binom{|S_1|}{k} \binom{|S_2|}{N_I - k}}{\binom{|S_1|+|S_2|}{N_I}}.
  \end{align}

  Write $m = |S_1| + |S_2|$, $|S_1| = \alpha m$ and   $N_I = \beta m$.
  By assumption $\alpha, \beta \in [\deltanew , 1-\deltanew] $ and $m \geq 2\deltanew (n+1)$. 
  Let $Z = \frac{N^1_I}{\alpha \beta m} -1$. 
  Under $\psap_{n+1} \left( \cdot \mcond \Gamma \in \vaso(\ga) \right)$, 
  $Z$ is a random variable of mean $0$,
  such that $\alpha \beta (1 + Z) m \in \ZZ \cap [0, \min \{ |S_1|,T_I \} ]$.
  
  First, we investigate the case where $Z$ is close to its mean.
  By means of a computation which uses Stirling's formula and \eqref{e.Tdistrib},  we find that
  \begin{align}\label{e.zdistrib}
 \psap_{n+1}\left(Z = z \mcond  \Gamma \in \vaso(\ga)  \right) 
    = \big( 1 + o(1) \big)
    \frac{ \exp \left(- \frac{\alpha \beta}{2(1-\alpha)(1-\beta)} m z^2 \right)}
    {\sqrt{2 \pi \alpha \beta (1 - \alpha)(1 - \beta) m}} \, ,     
 \end{align}
  where $o(1)$ designates a quantity tending to $0$ as $n$ tends to infinity, 
  uniformly in the acceptable choices of $\gamma$, $S_1$, $S_2$ and $z$, 
  with $|z| \leq \frac{2 n^{1/2}(\log n)^{1/2 +\eps}}{\alpha \beta m}$. We have obtained Lemma~\ref{l.balanced}(1).

  We now turn to the deviations of $Z$ from its mean. 
  From \eqref{e.Tdistrib}, one can easily derive that
  $\psap_{n+1}\left(Z = z \mcond \Gamma \in \vaso(\ga) \right)$ 
  is unimodal in $z$ with maximum at the value closest to $0$ that $Z$ may take.
  (We remind the reader that $Z$ takes values in $\frac{1}{\alpha \beta  m } \Z - 1$, 
  which contains $0$ only if ${\alpha \beta  m } \in \Z$.)  
The asymptotic equality
\eqref{e.zdistrib} thus implies the existence of constants $\mathrm{c}_0, \mathrm{c}_1 > 0$ depending only on $\deltanew$ such that, for $|z| \geq \frac{n^{1/2}(\log n)^{g}}{\alpha \beta m}$ and $n$ large enough, 
  \begin{align}\label{e.balanced0}
    \psap_{n+1} \Big( Z = z  \, \Big\vert \, \Gamma \in \vaso(\ga) \Big) \leq
  \mathrm{c}_1^{-1}  n^{-1/2}    \exp\left\{ - \mathrm{c}_0 ( \log n)^{2g}\right\}   \, ;
  \end{align}
  while for given $\eps > 0$, $|z| \geq \frac{n^{1/2}(\log n)^{1/2 + \eps}}{\alpha \beta m}$ and $n$ large enough,
  \begin{align}\label{e.balanced1}
    \psap_{n+1} \Big( Z = z \, \Big\vert \, \Gamma \in \vaso(\ga) \Big) \leq
   \mathrm{c}_1^{-1}  n^{-1/2}    \exp\left\{ - \mathrm{c}_0 ( \log n)^{1 + 2\eps}\right\}  \, \leq n^{- 2} \, .
  \end{align}

  Since $T^1_I$ takes no more than $n+1$ values, \eqref{e.balanced0} and~(\ref{e.balanced1}) imply Lemma~\ref{l.balanced}(2).\qed

\subsection{Mixing patterns by a random resampling}

Consider a random resampling experiment whose law we will denote by~$\pexper$. 
First an input polygon $\gamin$ is sampled according to the law $\psap_{n+1}$.
Then the contents of the slots in $S_1 \cup S_2$ are forgotten and independently resampled to form an output polygon $\gamout \in \sap_{n + 1}$. That is, given $\gamin$, $\gamout$ is chosen uniformly among the set of polygons $\ga \in \sap_{n + 1}$ for which $\ga \in \vaso(\gamin)$. Explicitly, if there are $j$ type $II$ patterns among $k$ slots in $S_1 \cup S_2$ in $\gamin$ (so that $k \geq j$), the polygon $\gamout$ is formed by choosing uniformly at random a subset of cardinality $j$ of these $k$ slots and inserting type $II$ patterns into the chosen slots.

Note the crucial property that $\gamout$ under $\pexper$ has the law $\psap_{n + 1}$: the resampling experiment holds the length~$n + 1$ random polygon at equilibrium. We mention in passing that a basic consequence of this resampling is a delocalization of the walk midpoint.
\begin{proposition}{\cite[Proposition 1.3]{ontheprob}}
Let $d \geq 2$. There exists $C > 0$ such that, for $m \in \N$,
$$
 \sup_{x \in \Z^d} \psaw_{m} \Big( \Ga_{\lfloor m/2 \rfloor} = x \Big) \leq C m^{-1/2} \, .
$$
\end{proposition}
It may be instructive to consider how to proof this result using the resampling experiment (or in fact a similar one involving walks rather than polygons) and Lemma~\ref{l.kptconseq}; a proof using such an approach is given in \cite[Section 3.2]{ontheprob}.

Any element $\phi \in \vaso(\ga)$ 
begins by tracing a journey over the region where slots in $S_1$ may appear, from the origin to $\gaemp_{(n+1)/10}$; it then follows its {\em middle section}, the trajectory of $\ga$ from  $\gaemp_{(n+1)/10}$ until  $\gaemp_{l' - (n+1)/10}$ (where $l' = \ellempty$); and it ends by moving from this vertex back to the origin, through the territory of slots in~$S_2$. 
Note that, in traversing the middle section, $\phi$ is exactly following a sub-walk of $\ga$, because no pattern changes have been made to this part of~$\ga$.  The timing of the middle section of this trajectory is advanced or retarded according to how many type~$II$ patterns are placed in the slots  in $S_1$.
Each extra such pattern retards the schedule by two units.
When $\phi$ has the {\em minimum} possible number $m : = \max \{ 0 , T_{II}(\ga)  - \vert S_2 \vert  \}$ of type $II$ patterns in the slots of $S_1$, the middle section is traversed by $\phi$ as early as possible,  the journey taking place during $\big[(n+1)/10 + 2m, n - 2 \big( T_{II}(\ga) - m \big) - (n+1)/10 \big]$. When $\phi$ has the maximum possible number $M : = \min \{ \vert S_1 \vert , T_{II}(\ga) \}$ of type~$II$ patterns in the slots of $S_1$, this traversal occurs as late as possible,  during $\big[ (n+1)/10 + 2 M , (n+1)  - (n+1)/10 - 2 \big( T_{II}(\ga) - M \big) \big]$. Since
$M \leq \vert S_1 \vert \leq (n+1)/13$ and $T_{II}(\ga) - m  \leq \vert S_2 \vert \leq (n+1)/13$, $\phi_j$ necessarily lies in the middle section whenever $j \in [(n+1)/10 + 2(n+1)/13,n+1 - (n+1)/10 -2(n+1)/13]$. Since the snake method parameter $\ell$ has been set to belong to the interval $[ \lceil n/4 \rceil , \lfloor 3n/4 \rfloor]$, we see that $\phi_j$ always lies in the middle section whenever $j \in \big[ \ell - n^{1/2}, \ell + n^{1/2} \big]$.

Taking $\ga \in \sap_{n + 1}$ and conditioning $\pexper$ on $\gamin = \ga$, note that the mean number of type $I$ patterns that end up in the slots in $S_1$ under $\gamout$ is given by $T_I(\ga) \cdot \tfrac{\vert S_1(\ga) \vert}{\vert S_1(\ga) \vert + \vert S_2(\ga) \vert}$, because this expression is the product of the number of type~$I$ patterns that are redistributed and the proportion of the available slots that lie in~$S_1$.

Consider now a polygon $\phi \in \vaso(\ga)$ that achieves as closely as possible the mean value for the number of type $I$ patterns among the slots in $S_1$: that is, $T_I^1(\phi)$ equals $\lfloor T_I(\ga) \cdot \tfrac{\vert S_1(\ga) \vert}{\vert S_1(\ga) \vert + \vert S_2(\ga) \vert} \rfloor$. As we have noted, $\phi_{\ell}$ is always reached during the middle section of $\phi$'s three-stage journey. Define  the {\em middle index}  $\lmid = \lmid(\ga)$ so that $\phi_{\ell} = \ga_{\lmid}$.   
Note that, given $\ga$ and this value of  $T_I^1(\phi)$, the value of this index is independent of the choice of $\phi$.

\subsection{Snakes of walks with high closing probability are typical}

Recall that the index parameter $\ell$ (and also $n$)
were fixed in Subsection~\ref{s.indexset}.
Moving towards the proof of Proposition~\ref{p.snakeg}, we take $\ga \in \sap_{n + 1}$ and define $\notick(\ga)$ to be the set  
\begin{eqnarray*}
     & & \Big\{ j \in \big(\ell - 2\N\big) \, \cap \, \big[ \lmid(\ga) - 2 n^{1/2} (\log n)^{1/4} \, , \,  \lmid(\ga) + 2 n^{1/2} (\log n)^{1/4} \big]: \\
   & & \qquad \qquad \qquad \qquad  \qquad \qquad  \qquad 
    \ga \, \, \textrm{is not $(1/2 - 2\eps,n,\ell)$-charming at $j$} \Big\} \, .
\end{eqnarray*}
Henceforth in this proof, charming will mean $(1/2 - 2\eps,n,\ell)$-charming. 
(The parity constraint that $\ell - j$ is even is applied above because the walks of length $j + n - \ell$ considered in the definition of charming at index~$j$ must be of odd length if they are to close.)

 We have seen that 
$\gamout_{\ell}$ is visited by $\gamout$ during its middle section, when it is traversing a subpath of $\gamin$ unchanged by pattern mixing. For this reason, we may define a random variable $L$ under $\pexper$ by setting  
$\gamout_{\ell} = \gamin_L$. 

We now state a key property of the resampling procedure. 
 \begin{lemma}\label{l.keyprop}
 The events $\big\{ \textrm{$\gamout$ is charming at $\ell$} \big\}$ and $\big\{ \textrm{$\gamin$ is charming at $L$} \big\}$ coincide.
\end{lemma}
\noindent{\bf Proof.} Note that the shells of $\vas\big( \gamin_{[0,L]} \big)$ and $\vas \big( \gamout_{[0,\ell]} \big)$ coincide, because   
$\gamout_{[0,\ell]}$ may be obtained from  $\gamin_{[0,L]}$ by modifying the $I/II$-status of some of its slots (these being certain slots in $S_1$). Thus,
Lemma~\ref{l.sharedavoidance} implies the statement. \qed

\begin{lemma}\label{l.jbound}
$$
\pexper \Big(   \big\vert \notick(\gamin) \big\vert \geq n^{1/2 - \eps/6}  \Big) \leq n^{-\eps/6} \, .
$$
\end{lemma}
\noindent{\bf Proof.}
Choosing $\deltanew > 0$ small enough and abbreviating $\mathcal{G} = \mathcal{G}_{n+1,\deltanew}$, Lemma~\ref{l.balanced}(1) implies that, 
for each $\ga \in \mathcal{G}$ and $k \in \big[  - n^{1/2} (\log n)^{1/4}, n^{1/2} (\log n)^{1/4} \big]$,
$$
 \pexper \Big(  L  = \lmid(\gamin) +  2k    \,  \Big\vert \,  \gamin =  \ga \Big) \geq  n^{-1/2} \exp \big\{ - c (\log n)^{1/2} \big\} \, .
$$ 
Thus, again taking any $\gamma \in \mathcal{G}$,
\begin{eqnarray*}
 & & \pexper \Big( \gamout \, \, \textrm{is not charming at $\ell$} \, \Big\vert \, \gamin  = \ga \Big) \\
 & \geq & \sum_{k =   - n^{1/2} (\log n)^{1/4}}^{n^{1/2} (\log n)^{1/4}}  \pexper \Big( \gamout \, \, \textrm{is not charming at $\ell$} \, , \, L =  \lmid(\gamin) + 2k \,  \Big\vert \, \gamin  = \ga \Big) \\
 & \geq & \sum_{k =   - n^{1/2} (\log n)^{1/4}}^{n^{1/2} (\log n)^{1/4}}  \pexper \Big(  \ga \, \, \textrm{is not charming at $\lmid(\gamin) + 2k$} \,   \, , \, L =  \lmid(\gamin) + 2k  \,  \Big\vert \, \gamin  = \ga \Big) \\
 & \geq &   n^{-1/2} \exp \big\{ - c (\log n)^{1/2} \big\} \sum_{k =   - n^{1/2} (\log n)^{1/4}}^{n^{1/2} (\log n)^{1/4}}  1\!\! 1_{    \ga \, \, \textrm{is not charming at $\lmid(\ga) + 2k$}}   \\
 & \geq &   n^{-1/2} \exp \big\{ - c (\log n)^{1/2} \big\} \cdot    \vert \notick(\ga) \vert  \, ,
\end{eqnarray*}
where the second inequality made use of Lemma~\ref{l.keyprop}.

Averaging over such $\gamma$, we find that
\begin{eqnarray*}
 & & \pexper \Big(  \gamout \, \, \textrm{is not charming at $\ell$} \, \Big\vert \, \gamin  \in \mathcal{G} \Big) \\
& \geq & c n^{-1/2}  \exp \big\{ - c (\log n)^{1/2} \big\}  \cdot \eexper \Big[ \,  \big\vert \notick(\gamin) \big\vert \, \Big\vert \, \gamin  \in \mathcal{G} \Big] \, ,
\end{eqnarray*}
where $\eexper$ denotes the expectation associated with the law $\pexper$.

Note that
\begin{eqnarray*}
 & & \pexper \Big( \gamout \, \, \textrm{is not charming at $\ell$} \, \Big\vert \, \gamin  \in \mathcal{G} \Big) \\
 & \leq &  2 \, \psap_{n+1} \Big(  \Gamma \, \, \textrm{is not charming at $\ell$} \Big) 
  =  2 \, \psap_{n+1} \Big( \Gamma_{[0,\ell]} \not\in \hPhi_{\ell,n}^{1/2 - 2\eps}  \Big) \leq 2 n^{-\eps} \, ,
\end{eqnarray*}
where in the first inequality we use that  $\gamout$ under $\pexper$ has the law $\psap_{n+1}$, and then apply Lemma~\ref{l.kptconseq}, to find that $\pexper \big( \gamin \in \mathcal{G} \big) = \psap_{n+1} \big( \Gamma \in \mathcal{G} \big) \geq 1 - e^{-cn} \geq 1/2$. The final inequality above used Lemma~\ref{l.polyclose}.

Thus,
$$
\eexper \Big[ \,  \big\vert \notick(\gamin) \big\vert \, \Big\vert \, \gamin  \in \mathcal{G} \Big] \leq  n^{1/2 - \eps/2} \, .
$$

We find that $\eexper \vert \notick \vert$ is at most
\begin{eqnarray*}
  &  & 
\eexper \Big[ \,  \big\vert \notick(\gamin) \big\vert \, \Big\vert \, \gamin  \in \mathcal{G} \Big] \, + \, \Big( 2 n^{1/2} \big( \log n \big)^{1/4} + 1 \Big)  \pexper \Big( \gamin \not\in \mathcal{G} \Big) \\
   & \leq & n^{1/2 - \eps/2} \, + \, \Big( 2 n^{1/2} \big( \log n \big)^{1/4} + 1 \Big) \, e^{-cn}  \leq n^{1/2 - \eps/3} \, . 
\end{eqnarray*}

Applying Markov's inequality yields Lemma~\ref{l.jbound}. \qed

\medskip

\noindent{\bf Proof of Proposition~\ref{p.snakeg}.}
Lemma~\ref{l.balanced}(2) implies that 
$$
\psap_{n+1} \Big( \ell \in \big[ \lmid(\Ga) -  n^{1/2} (\log n)^{1/4} ,  \lmid(\Ga) +  n^{1/2} (\log n)^{1/4} \big] \, \Big\vert \,  \Ga \in \mathcal{G}_{n,\deltanew} \Big) 
$$
is at least $1 - \exp \big\{ - c (\log n)^{1/2} \big\}$. Note that the interval centred on $\lmid(\Ga)$ considered here is shorter than its counterpart in the definition of $\notick(\ga)$ for $\ga \in \sap_{n+1}$.
Applying Lemmas~\ref{l.kptconseq} and~\ref{l.jbound}, we find that
\begin{eqnarray*}
 & &  \psap_{n+1} \Big( \# \Big\{  j \in \big[ \ell - n^{1/2} (\log n)^{1/4} ,  \ell + n^{1/2} (\log n)^{1/4} \big]: \\
   & & \qquad \qquad   \textrm{$\Gamma$ is not charming at $j$} \Big\} \geq n^{1/2 - \eps/6}   \Big) 
 \leq   n^{-\eps/6} +  e^{- c (\log n)^{1/2}}    + e^{-cn} \, .
\end{eqnarray*}

When the complementary event occurs, $\Gamma$ is charming for at least one-quarter of the indices in  $\big[ \ell - n^{1/2}, \ell  \big]$, so that $\Gamma_{[0,\ell]}$ is an element of $\mathsf{CS}_{1/2,0}^{1/2 - 2\eps,\ell,n}$. (We write one-quarter rather than one-half here, because one-half of such indices are inadmissible due to their having the wrong parity.) \qed

\bibliographystyle{plain}

\bibliography{saw}

\begin{thebibliography}{10}

\bibitem{BBSlog}
R.~Bauerschmidt, D.~Brydges, and G.~Slade.
\newblock Logarithmic correction for the susceptibility of the 4-dimensional
  weakly self-avoiding walk: a renormalisation group analysis.
\newblock arXiv:1403.7422. Commun. Math. Phys., to appear, 2014.

\bibitem{BDGS11}
R.~Bauerschmidt, H.~Duminil-Copin, J.~Goodman, and G.~Slade.
\newblock Lectures on self-avoiding walks.
\newblock In D.~Ellwood, C.~Newman, V.~Sidoravicius, and W.~Werner, editors,
  {\em Lecture notes, in Probability and Statistical Physics in Two and More
  Dimensions}. CMI/AMS -- Clay Mathematics Institute Proceedings, 2011.

\bibitem{BMB09}
Mireille Bousquet-M{\'e}lou and Richard Brak.
\newblock Exactly solved models.
\newblock In {\em Polygons, polyominoes and polycubes}, volume 775 of {\em
  Lecture Notes in Phys.}, pages 43--78. Springer, Dordrecht, 2009.

\bibitem{ontheprob}
Hugo Duminil-Copin, Ioan Manolescu, Alexander Glazman, and Alan Hammond.
\newblock On the probability that self-avoiding walks ends at a given point.
\newblock arXiv:1305.1257. Ann. Probab., to appear, 2013.

\bibitem{Dup89}
B.~Duplantier.
\newblock Fractals in two dimensions and conformal invariance.
\newblock {\em Phys. D}, 38(1-3):71--87, 1989.
\newblock Fractals in physics (Vence, 1989).

\bibitem{Dup90}
B.~Duplantier.
\newblock Renormalization and conformal invariance for polymers.
\newblock In {\em Fundamental problems in statistical mechanics {VII}
  ({A}ltenberg, 1989)}, pages 171--223. North-Holland, Amsterdam, 1990.

\bibitem{Flory}
P.~Flory.
\newblock {\em Principles of Polymer Chemistry}.
\newblock Cornell University Press, 1953.

\bibitem{hammersleywelsh}
J.~M. Hammersley and D.~J.~A. Welsh.
\newblock Further results on the rate of convergence to the connective constant
  of the hypercubical lattice.
\newblock {\em Quart. J. Math. Oxford Ser. (2)}, 13:108--110, 1962.

\bibitem{CJC}
Alan Hammond.
\newblock On self-avoiding polygons and walks: counting, joining and closing.
\newblock arXiv:1504.05286, 2017.

\bibitem{snakemethodpolygon}
Alan Hammond.
\newblock On self-avoiding polygons and walks: the snake method via polygon
  joining.
\newblock \url{math.berkeley.edu/~alanmh/papers/snakemethodpolygon.pdf}, 2017.

\bibitem{ThetaBound}
Alan Hammond.
\newblock An upper bound on the number of self-avoiding polygons via joining.
\newblock \url{math.berkeley.edu/~alanmh/papers/ThetaBound.pdf}, 2017.

\bibitem{HS92}
T.~Hara and G.~Slade.
\newblock Self-avoiding walk in five or more dimensions. {I}. {T}he critical
  behaviour.
\newblock {\em Comm. Math. Phys.}, 147(1):101--136, 1992.

\bibitem{HS92b}
Takashi Hara and Gordon Slade.
\newblock The lace expansion for self-avoiding walk in five or more dimensions.
\newblock {\em Rev. Math. Phys.}, 4(2):235--327, 1992.

\bibitem{JROSTW}
E.~J. Jance~van Rensburg, E.~Orlandini, D.~W. Sumners, M.~C. Tesi, and S.~G.
  Whittington.
\newblock The writhe of a self-avoiding polygon.
\newblock {\em Journal of Physics A: Mathematical and General}, 26(19):L981 --
  L986, 1993.

\bibitem{kestenone}
H.~Kesten.
\newblock On the number of self-avoiding walks.
\newblock {\em J. Mathematical Phys.}, 4:960--969, 1963.

\bibitem{Lawler13}
G.~Lawler.
\newblock Random walk problems motivated by statistical physics.
\newblock http://www.math.uchicago.edu/~lawler/russia.pdf, 2013.

\bibitem{LSW}
G.~F. Lawler, O.~Schramm, and W.~Werner.
\newblock On the scaling limit of planar self-avoiding walk.
\newblock In {\em Fractal geometry and applications: a jubilee of {B}eno\^\i t
  {M}andelbrot, {P}art 2}, volume~72 of {\em Proc. Sympos. Pure Math.}, pages
  339--364. Amer. Math. Soc., Providence, RI, 2004.

\bibitem{Madras14}
N.~Madras.
\newblock A lower bound for the end-to-end distance of the self-avoiding walk.
\newblock {\em Canad. Math. Bull.}, 57(1):113--118, 2014.

\bibitem{MS93}
N.~Madras and G.~Slade.
\newblock {\em The self-avoiding walk}.
\newblock Probability and its Applications. Birkh\"auser Boston Inc., Boston,
  MA, 1993.

\bibitem{Nie82}
B.~Nienhuis.
\newblock Exact critical point and critical exponents of ${O}(n)$ models in two
  dimensions.
\newblock {\em Phys. Rev. Lett.}, 49:1062--1065, 1982.

\bibitem{Nie84}
B.~Nienhuis.
\newblock Coulomb gas description of {2D} critical behaviour.
\newblock {\em J. Statist. Phys.}, 34:731--761, 1984.

\bibitem{Orr47}
W.J.C. Orr.
\newblock Statistical treatment of polymer solutions at infinite dilution.
\newblock {\em Transactions of the Faraday Society}, 43:12--27, 1947.

\end{thebibliography}

\end{document}